# Short Note: Every Large Set of Integers Contains a Three Term Arithmetic Progression

GÁBOR KORVIN[1]

**Abstract**: I show that a trivial modification of a standard proof of the *Roth's Theorem* on triples in arithmetic progression would lead to the following *Theorem*: If $1 \leq n_1 < n_2 < \cdots$ is a "large set" that is if $\sum \frac{1}{n_i} = \infty$ then the sequence $\{n_1, n_2, \cdots\}$ contains an arithmetic progression of length three.

In 1936 Erdös and Turán conjectured [1] that if an infinite set of increasing positive integers $A = \{1 \leq n_1 < n_2 < \cdots\}$ has *positive upper density* then it contains arbitrary long arithmetic progressions (AP's). In 1952 Roth [2] proved, using the *circle method* of Hardy and Littlewood, that under this condition *A* always contains at least a 3-term AP. In 1962 Szemerédi [3] showed that *A* must contain at least a 4-term AP, and in 1974 [4] he succeeded to prove the full Erdös-Turán conjecture. Szemerédi's proof is combinatorial. A few years later Furstenberg [5, 6] found an independent proof of the Erdös-Turán conjecture, based on *Ergodic Theory*.

Erdös also posed [7, 8] the following much stronger problem (for the solution of which he originally offered 3,000 dollars): *Is it true that if $1 \leq n_1 < n_2 < \cdots$ and $\sum \frac{1}{n_i} = \infty$ then for any $n = 1, 2, \cdots$ the sequence $\{n_1, n_2, \cdots\}$ contains an arithmetic*

---

[1] Box 1157, Earth Sciences Department and Reservoir Characterization Research Group, King Fahd University, Dhahran 31261, Saudi Arabia. Email: gabor@kfupm.edu.sa

*progression of length n?* The conjecture would imply the recently proven theorem of Green and Tao [9] according to which the primes contain arithmetic progressions of an arbitrarily length *n*. (The longest progression found thus far was reported in January, 2007, by Jaroslaw Wrobleski, it consists of 24 primes [10] ).

Sets of integers $1 \leq n_1 < n_2 < \cdots$ with $\sum \frac{1}{n_i} = \infty$ are called *large sets*. In the light of the *Erdös conjecture*, it is surprising that the "largeness" of the set is not a necessary condition for its containing arbitrary long AP's. Indeed, it is easy to construct examples for *"small sets"* (i.e. sequences of integers $1 \leq n_1 < n_2 < \cdots$ such that $\sum \frac{1}{n_i} < \infty$) which contain arbitrarily long AP's. As an example, consider the sequence

$$S = \{1,10,11,\underbrace{100,101,102},\underbrace{1000,1001,1002,1003},\cdots,\underbrace{10^n,10^n+1,10^n+2,\cdots,10^n+n},10^{n+1},\cdots\}$$

For every $n = 1, 2, \cdots$ the sequence *S* contains AP's of length *n*. On the other hand, the series $\frac{1}{1} + \frac{1}{10} + \frac{1}{11} + \frac{1}{100} + \frac{1}{101} + \frac{1}{102} + \frac{1}{1000} + \cdots \frac{1}{1003} + \cdots$ is convergent. Indeed, one has

$$\frac{1}{1} + \frac{1}{10} + \frac{1}{11} + \frac{1}{100} + \frac{1}{101} + \frac{1}{102} + \frac{1}{1000} + \cdots \frac{1}{1003} + \cdots < \frac{1}{10^0} + \frac{2}{10^1} + \frac{3}{10^2} + \frac{4}{10^3} + \cdots$$

$$= 1 + 2x + 3x^2 + 4x^3 + \cdots \Big|_{x=0.1} = \frac{d}{dx}(x + x^2 + x^3 + \cdots)\Big|_{x=0.1} = \frac{d}{dx}\frac{x}{1-x}\Big|_{x=0.1} = \frac{1}{0.9^2}.$$

In this *Short Note* I call attention to the fact that a trivial modification of Newman's ([11], Chapter 4) elegant proof of the *Roth's Theorem* [2] would lead to the following

*Theorem*: If $1 \leq n_1 < n_2 < \cdots$ and $\sum \frac{1}{n_i} = \infty$ then the sequence $\{n_1, n_2, \cdots\}$ contains a nontrivial arithmetic progression of length 3.

We shall need the following simple *Lemma* [12-14]

*Lemma* 1: Let $\{m_1 < m_2 < \cdots\}$ be a sequence of integers with counting function f(n). Then

$$\sum_{n=1}^{\infty} \frac{1}{m_n} < \infty \; iff \; \sum_{n=1}^{\infty} \frac{f(n)}{n^2} < \infty. \qquad (1)$$

As Newman's book is easily available, I only recall here his notations and the main ideas he used in proving *Roth's Theorem*. Newman introduces the concept "*affine property*" P of finite sets of integers: A property P is *affine* if it satisfies the following two conditions:

C1. For each fixed pair of integers $\alpha \neq 0$ and $\beta$ the set $\{a_n\}$ has P if and only if $\{\alpha a_n + \beta\}$ has P;

C2. If a set has P, then all of its subsets have P.

(Examples are the trivial property $P_0$ of just being a set, or $P_A$ of not containing any arithmetic progressions of at least three distinct terms). Take now any fixed affine property P, and assume that the largest possible subset with property P, $S(n;P) \subseteq \{1,2,\cdots,n\}$ has $f(n;P) = |S(n;P)|$ elements. By Conditions C1 and C2, for any P the function $f(n;P)$ is *subadditive*, $f(m+n;P) \leq f(m;P) + f(n:P)$ that is (by Fekete's *Lemma*, [15]) $\lim_{n \to \infty} \frac{f(n,P)}{n} = C_P$ exists, and of course $0 \leq C_P \leq 1$.

Denote by $A(n;P)$ the number of arithmetic progressions that can be constructed from the elements of $S(n;P)$ where order counts and equal terms are allowed.

Newman proves Roth's Theorem through the following two *Lemmas*:

*Lemma 2 (Newman):* $\sum_{a \in S(n;P)} z^a = C_P \sum_{k \leq n} z^k + o(n)$ *uniformly on* $|z|=1$. (2)

*Lemma 3. (Newman):* $A(n;P) = \dfrac{C_P^3}{2} n^2 + o(n^2)$ (3)

Based on these two *Lemmas*, Newman states and proves Roth's Theorem as follows:

*Theorem (Roth):* $C_{P_A} = 0$.

Newman's proof: The only possible three term arithmetic arithmetic progressions formed from the elements of $S(n; P_A)$ are those whose three terms are equal, that is $A(n; P_A) \leq n$, that is by (Eq. 3)

$$\dfrac{C_{P_A}^3}{2} n^2 + o(n^2) \leq n,  \quad (4)$$

and letting $n \to \infty$ we get $C_{P_A} = 0$.

Make the following modification of the proof. Define $C_P(n) := \dfrac{f(n;P)}{n}$. It is easy to check that Newman's proofs of *Lemmas* 2 and 3 can be repeated without any other change than using $C_P(n)$ instead of $C_P$, so that we get – instead of *Lemmas* 2 and 3 - the following modified *Lemmas*:

*Lemma 2*$^*$: $\sum_{a \in S(n;P)} z^a = C_P(n) \sum_{k \leq n} z^k + o(n)$ *uniformly on* $|z|=1$. (5)

*Lemma 3*$^*$: $A(n;P) = \dfrac{C_P^3(n)}{2} n^2 + o(n^2)$ (6)

Now, we are ready to prove

*Theorem*: If $1 \leq n_1 < n_2 < \cdots$ and $\sum \dfrac{1}{n_i} = \infty$ then the sequence $\{n_1, n_2, \cdots\}$ contains a nontrivial arithmetic progression of length 3.

Proof : Assume that $n_1 < n_2 < \cdots$ has the property $P_A$ of not containing any arithmetic progressions of at least three distinct terms. As above (in the original Newman's proof of Roth's Theorem), we get from *Lemmas 2\** and *3\** (writing $\frac{f(n;P_A)}{n}$ instead of $C_{P_A}(n)$): $\left[\frac{f(n,P_A)}{n}\right]^3 \frac{n^2}{2} + o(n^2) \leq n$ which implies for $n$ greater than some sufficiently large $n_0$ that $\left[\frac{f(n;P_A)}{n}\right]^3 \frac{n^2}{2} \leq n$. Easy algebra gives that in this case $\frac{f(n;P_A)}{n^2} \leq \frac{\sqrt[3]{2}}{n^{4/3}}$. Consequently, $\sum_{n=1}^{\infty} \frac{f(n;P_A)}{n^2} \leq \sum_{n=1}^{n_0} \frac{f(n;P_A)}{n^2} + \sqrt[3]{2} \sum_{n=n_0+1}^{\infty} \frac{1}{n^{4/3}} < \infty$ because (see e.g. [16], p. 115) $\sum_{n=1}^{\infty} \frac{1}{n^{\alpha}} < \infty$ for $\alpha > 1$. By *Lemma 1* we get $\sum_{i=1}^{\infty} \frac{1}{n_i} < \infty$. The contradiction shows that the *sequence $\{n_1, n_2, \cdots\}$ cannot have the property $P_A$, but it contains an arithmetic progression of length 3.*

*Acknowledgment* I gratefully acknowledge the peaceful, creative atmosphere of the King Fahd University of Petroleum and Minerals, that has made it possible for me to leave aside for a brief time my more practical tasks in rock physics and make a small *excursus* to *Combinatorial Number Theory*.